\documentclass{article}

\usepackage{amssymb,amsthm,amsmath,amsfonts,amstext}
\usepackage{url}
\usepackage{latexsym}
\begin{document}

\title{A curious proof of Fermat's little theorem}
\author{Giedrius Alkauskas}
\date{}
\maketitle

\indent Fermat's little theorem states that for $p$ prime and $a\in\mathbb{Z}$,
$p$ divides $a^{p}-a$. This result is of huge importance in elementary and
algebraic number theory. For instance, with its help we obtain the so-called
Frobenius automorphism of a finite field $\mathbb{F}_{p^{n}}$ over
$\mathbb{F}_{p}$.\\
\indent This theorem has many interesting and sometimes unexpected proofs. One
classical proof is based upon properties of binomial coefficients. In fact, $
(d+1)^{p}-d^{p}-1= \sum_{i=1}^{p-1}\binom{p}{i}d^{i}. $ Since
$\binom{p}{i}=\frac{p!}{i!(p-i)!}$ is divisible by $p$ for $1\leq i\leq p-1$,
then $(d+1)^{p}-d^{p}-1$ is divisible by $p$. Summing this over
$d=1,2,...,a-1$, we obtain the desired result. Another classical proof is based
upon Lagrange's theorem, which states that the order of an element of a finite
group divides the group order. Applying this theorem to the multiplicative
group of a finite field $\mathbb{F}_{p}$ we obtain the result immediately.
Several other proofs can be found at \cite{wiki}.
 Nevertheless, in all of these proofs one or another analogue of the Euclidean algorithm
 (hence arithmetic) is being used.\\
\indent In this short note we present a curious proof which was found as a side
result of another, unrelated problem (which is the case, maybe, with many such
``curious" proofs). Surprisingly, arithmetic, algebra,
and the properties of binomial coefficients do not manifest at all.\\

Let $f(x)=1-x-dx^{2}+\sum_{k\geq3}a_{k}x^{k}$ be any formal power series in
$\mathbb{Q}$, with coefficients in $\mathbb{Z}$. It is well known that this
series can be represented in a unique way as a formal product of the following
form:
\[
f(x)=\prod_{k\geq1}(1-m_{k}x^{k}),
\]
where the coefficients $m_{k}$ are integers. This result can be found in
\cite{kobliz}, but the proof is simple and straightforward. In fact, for $k=1$
and $k=2$ we have a unique choice $m_{1}=1$ and $m_{2}=d$. Suppose $N\geq 3$
and we have already chosen $m_{k}$ for $k\leq N-1$. Then
$\prod_{k=1}^{N-1}(1-m_{k}x^{k})=1-x-dx^{2}+\sum_{k=3}^{N-1}a_{k}x^{k}+Cx^{N}+\text{``higher
terms"}$, where $C$ is a certain integer which depends only on $m_{k}$ for
$1\leq k\leq N-1$. Therefore, the unique choice for $m_{N}$ is $m_{N}=C-a_{N}$.
In a similar fashion, since
$\frac{1}{f(x)}=1+x+(d+1)x^{2}+\sum_{k\geq3}b_{k}x^{k}$ is also a formal
integer power series, it can be represented in a unique way as a product
\[
\frac{1}{f(x)}=(1+x)(1+(d+1)x^{2})\prod_{k\geq3}(1-n_{k}x^{k}),
\]
where $n_{k}$ are integers as well, $n_{1}=-1$, and $n_{2}=-(d+1)$.\\
\indent Recall that the logarithmic derivative of a power series $g(x)$,
denoted by $(\ln g(x))'$, is defined to be the power series $g'(x)/g(x)$. It is
not hard to prove that for any two formal power series $g(x)$ and $h(x)$, $(\ln
g(x)\cdot h(x))'=(\ln g(x))'+(\ln h(x))'$. Indeed, this property reduces to the
Leibniz rule
\begin{eqnarray*}
(g(x)\cdot h(x))'=g'(x)h(x)+g(x)h'(x).
\end{eqnarray*}
This is verified simply by comparing the corresponding coefficients. Note also
that the binomial theorem is not used in the proof.\\
\indent Now take the formal logarithmic derivative of $f(x)$. We obtain:
\[
-x\Big{(}\ln f(x)\Big{)}'=\sum_{k\geq1}\frac{km_{k}x^{k}}{1-m_{k}x^{k}}=
\sum_{N\geq1}x^{N}\sum_{s|N}m_{N/s}^{s}\frac{N}{s}.
\]
In a similar fashion,
\[
-x\Big{(}\ln\frac{1}{f(x)}\Big{)}'=x(\ln f(x))'=\sum_{N\geq 1}x^{N}
\sum_{s|N}n_{N/s}^{s}\frac{N}{s}.
\]
Therefore, we have interesting identities among the terms of two infinite
sequences:
\begin{eqnarray}
\sum_{s|N}m_{N/s}^{s}\frac{N}{s}=- \sum_{s|N}n_{N/s}^{s}\frac{N}{s},\quad
N\in\mathbb{N}.\label{rec}
\end{eqnarray}
\indent We can easily prove by induction that this implies $m_{k}=-n_{k}$ for
odd $k$, but not for the terms with even indices! Thus, a consequence of this
reasoning is the fact that any infinite sequence of integers
$\{m_{k},k\in\mathbb{N}\}$ with $m_{1}=\pm 1$ has an ``inverse" sequence of
integers $\{n_{k},k\in\mathbb{N}\}$ with $n_{1}=\mp1$. Consequently, all such
sequences split into mutually inverse pairs. It is rather tempting to try to
express an inverse of a certain sequence for which the infinite product has a
rich mathematical content. For example, let us take $m_{k}=1$ for
$k\in\mathbb{N}$. Hence, we have a product
\begin{eqnarray*}
(x,x)_{\infty}=\prod\limits_{k=1}^{\infty}(1-x^{k}).
\end{eqnarray*}
It is well known that
$(x,x)_{\infty}^{-1}=\sum\limits_{n=0}^{\infty}p(n)x^{n}$, where $p(n)$ is
Ramanujan's partition function. Using the recurrence (\ref{rec}) we can compute
the sequence $\widetilde{n}_{k}=-n_{k}$. As mentioned, $\widetilde{n}_{k}=1$
for $k$ odd, and terms of this sequence with even indices begin with
\begin{eqnarray*}
2,4,0,14,-4,-8,-16,196,-54,-92,-184,144,-628,-1040,-2160,41102...
\end{eqnarray*}
Therefore,
\begin{eqnarray*}
\sum\limits_{n=0}^{\infty}p(n)x^{n}=\prod\limits_{k=1}^{\infty}(1+\widetilde{n}_{k}x^{k}).
\end{eqnarray*}

Let us return to our case. Recall that $m_{2}=d$ and $n_{2}=-(d+1)$. Hence,
when $N=2p$, where $p>2$ is a prime, (\ref{rec}) reads as:
\[
2p\cdot m_{2p}+p\cdot m_{p}^{2}+2d^{p}+1= -2p\cdot n_{2p}-p\cdot
n_{p}^{2}+2(d+1)^{p}-1.
\]
Thus, $p$ divides $(d+1)^{p}-d^{p}-1$. Summing this over
$d=1,2,...,a-1$, we finally obtain $p|a^{p}-a$. Quite unexpected!\\

Likewise, expand the following function into a formal infinite product:
\[
f(x)=1-x-\sum_{n=1}^{\infty}d^{n}x^{n+1}=\prod_{n=1}^{\infty}(1-a_{n}x^{n}).
\]
Since $f(x)=\frac{1-(d+1)x}{1-dx}$, after taking the logarithmic derivative, we
obtain:
\[
-x\Big{(}\ln
f(x)\Big{)}'=\sum\limits_{N=1}^{\infty}\Big{(}(d+1)^{N}-d^{N}\Big{)}x^{N}=
\sum_{N\geq1}x^{N}\sum_{s|N}a_{N/s}^{s}\frac{N}{s}.
\]
As a direct consequence, $a_{p}=\frac{(d+1)^{p}-d^{p}-1}{p}$, which implies
that $\frac{(d+1)^{p}-d^{p}-1}{p}$ is an integer. Possible variations on this
theme unexpectedly produce other congruences and identities. Recall that a
prime number $p$ is said to be a Wieferich prime if and only if
$2^{p-1}\equiv1\text{ (mod }p^{2})$. Examples are $p=1093$ and $p=3511$, with
no others in the range $p<4\cdot10^{12}$. In the last example with $d=1$, {\it
all} the numbers $a_{p}=\frac{2^{p}-2}{p}$ appear simultaneously in the
infinite product defining $\frac{1-2x}{1-x}$, and as the proof of the algorithm
used to expand a formal power series into an infinite product suggests,
strangely enough, the coefficients $a_{N}$ are defined inductively on $N$
without a distinction between prime and composite values of $N$. Possibly, more
profound research of this product could clarify our understanding of these
exceptional Wieferich primes.

\bigskip

\noindent\textit{Department of Mathematics and Informatics,
Vilnius University, Naugarduko 24, 03225 Vilnius, LITHUANIA\\
giedrius.alkauskas@maths.nottingham.ac.uk}

\end{document}